\documentclass[11]{article}

\newtheorem{Probl}{Problem}
\newtheorem{Prop}{Proposition}
\newtheorem{PropDef}{Proposition-Definition}
\newtheorem{Lem}{Lemma}
\newtheorem{Def}{Definition}

\begin{document}

\begin{center}
{\Huge\bf
Interval Computations and their Categorification
}\vspace{.5in}

{Nikolaj Glazunov}     \\
\smallskip

Glushkov Institute of Cybernetics NAS\\
40 Glushkov prospect\\
03187 Kiev GSP-680\\

glanm@d105.icyb.kiev.ua
\end{center}



\begin{abstract}
By the example of the proof of Minkowski's conjecture on critical determinant
we give a category theory framework for interval computation.
\end{abstract}

\section{Introduction}
The purpose of this paper is to describe interval computations
under the proof of Minkowski's conjecture on the critical
determinant of the region $ |x|^p + |y|^p < 1, \; p > 1 $
 in terms of (certain kind of) categories and functors.
Let $ D_p \subset {\bf R}^2 = (x,y), \ p > 1 $ be the 2-dimensional
region:
$$ |x|^p + |y|^p < 1 . $$
The well known Minkowski conjecture~\cite{Mi:DA,M:LP,D:NC,Co:MC,W:MC,Ma:AC}
about critical determinant of
the region $ D_p $ can be reformulated as the problem of minimization
on moduli space $\cal M $ of admissible lattices of the region
$ D_p$~\cite{G:RM}.
 The moduli space is defined by the equation
$$ \Delta(p,\sigma) = (\tau + \sigma)(1 + \tau^{p})^{-\frac{1}{p}}
  (1 + \sigma^p)^{-\frac{1}{p}}, \; \; \;  \; (1) $$
in the domain
 $$ D_{p}: \; \infty > p > 1, \; 1 \leq \sigma \leq \sigma_{p} =
 (2^p - 1)^{\frac{1}{p}}, $$
of the $ \{p,\sigma\} $ plane, where $\sigma$ is some real parameter;
$\;$ here $ \tau = \tau(p,\sigma) $ is the function uniquely
determined by the conditions
$$ A^{p} + B^{p} = 1, \; 0 \leq \tau \leq \tau_{p}, $$
where
$$ A = A(p,\sigma) = (1 + \tau^{p})^{-\frac{1}{p}} -
(1 + \sigma^p)^{-\frac{1}{p}},            \;
 B = B(p,\sigma) = \sigma(1 + \sigma^p)^{-\frac{1}{p}}       +
\tau(1 + \tau^{p})^{-\frac{1}{p}}, $$ $\tau_{p}$ is defined by the
equation $ 2(1 - \tau_{p})^{p} = 1 + \tau_{p}^{p}, \; 0 \leq
\tau_{p} < 1. $ \\ {\bf Minkowski's analytic conjecture:}
\\ { \it For any real $p$ and $\tau$ with conditions } $ p > 1, \
p \ne 2, \ 0 < \tau < \tau_p $ $$ \Delta(p,\sigma) >
min(\Delta(p,1),\Delta(p,\sigma_p)). \; \; (2) $$ A.V. Malyshev
and I have proved by interval computation \\ {\bf Theorem.} {\it
For all $p \geq 1.01, \; p \neq 2$ and for all $ 1 < \sigma <
\sigma_{p}$ with the exception of $ \{p,\sigma\} $ from the domain
$P_{2,\sqrt(3)} = \{2 \leq p \leq 2.000003, \; d \leq \sigma \leq
\sigma_{p}\} \; d = 1.7320503$ the inequality (2) takes place. }
\\
 Nondifferential interval method of computation of implicitly defined
functions is one of the ingredients of this proof. It can be
interpreted as a contracting map on the interval space. The map is
defined by some paths of the corresponding program. \\
  The category theory view on manifolds in mathematics and on programs
 in computer science as well as interval computations on Minkowski
moduli space led us to the necessity of the introduction of interval
manifolds, presheaves on them and functors from path categories of
programs to the interval categories. Using these and some other notions,
we  present a category theory framework for interval computation.
Proofs are omitted.
\section{Lattices, Admissible Lattices and Critical determinants}
   In this section we recall notions of lattices, admissible lattices
and critical determinants and some of their properties~\cite{C:GN}.

 Let $a_{1}, \ldots ,a_n $ be the independent points (a basis) of ${\bf R}^n$.
The set $\Lambda$  of points
$$ x = u_{1}a_{1} + u_{2}a_{2}+ \ldots +u_{n}a_{n}, \;
(u_{1}, \ldots ,u_{n} \; integers)  $$
is called a lattice. The system of points $a_{1}, \ldots ,a_n $
is called a {\it basis} of $\Lambda$. \\

   If $\Lambda$ is a lattice and $A$ is a basis of $\Lambda,$ then
$|det A|$ is called the {\it determinant} of $\Lambda$. It is denoted by
$d(\Lambda)$. \\

   Let $M$ be an arbitrary  set in ${\bf R}^n$, $O = (0,0) \in {\bf R}^n$.
A lattice
$\Lambda$ is called admissible for $M$, or $M-${\it admissible}, if
it has no points $\neq$ O in the interior of $M.$ It is called
{\it strictly admissible} for
$M$ if it does not contain a point $\neq$ O of $M.$     \\

   The {\it critical determinant} of a set $M$ is the quantity
$\Delta(M)$ given by
$$\Delta(M) = inf\{d(\Lambda): \Lambda \; strictly  \,
admissible \, for M\}  $$
with the understanding that $\Delta(M) = \infty$  if there are no
strictly admissible lattices. The set $M$ is said to be of the finite
or the infinity type according to whether $\Delta(M)$ is finite or
infinite. \\

  Of course, we may also consider the greatest lower bound of
$d(\Delta)$ on the collection of all $M-$admissible lattices.
We put
$$\Delta_{0}(M) = inf\{d(\Lambda): \Lambda \; admissible \, for M\}.$$
Then clearly, $\Delta_{0}(M) \leq \Delta(M)$. Here, for large classes
of sets, the equality sing hold. For the class of open sets, this is
trivial. As we consider open sets, we will denote the critical
determinants by $\Delta(M)$.

\section{Interval Cellular Covering}

  For any $n$ and any $j, \; 0 \leq j \leq n$, an $j-${\it dimensional
interval cell}, or $j-${\it I-cell}, in ${\bf R}^n$ is a subset $Ic$ of
${\bf R}^n$
such that (possibly, after permutation of variables) it has the form   \\

$ Ic = \{x \in {\bf R}^n: {\underline a}_{i},{\overline a}_{i}, r_{k} \in
{\bf R} :    \\
 \; {\underline a}_{i} \le x_{i} \le {\overline a}_{i}, 1 \le i \le j, \\
 \; x_{j + 1} = r_{1}, \cdots, x_{n} = r_{n-j} \} \; . $
Here $  {\underline a}_{i} \le {\overline a}_{i}. $  \\

If $j = n$ then we have an $n-$dimensional interval vector.
Let ${\cal P}$ be the hyperplane that contains $Ic.$ These is the
well known fact:
\begin{Lem}
 The dimension of $Ic$ is equal to the minimal dimension of hyperplanes
that contain $Ic.$
\end{Lem}
 Let ${\cal P}$ be the such hyperplane, $Int \; Ic$ the set of interior
points
of $Ic$ in ${\cal P}, \; Bd \; Ic = Ic \setminus Int \; Ic.$ For
$m-$dimensional I-cell $Ic$ let $d_{i}$ be an $(m-1)-$dimensional I-cell from
$Bd \; Ic.$ Then
$d_{i}$ is called an $(m-1)-$dimensional {\it face} of the I-cell $Ic.$
\begin{Def}
 Let $D$ be a bounded set in ${\bf R}^n.$ By interval cellular covering
$Cov$ we will understand any finite set of $n-$dimensional I-cells  such that
their union contains $D$ and adjacent I-cells are intersected by their
faces only. By $ {\mid Cov \mid}$ we will denote
the union of all I-cells from $Cov.$
\end{Def}
  Let $Cov$ be the interval covering. By its  {\em subdivision}
we will understand an interval covering $Cov^{`}$ such that
$\mid Cov \mid = \mid Cov^{`} \mid$ and each I-cell from $Cov^{`}$ is
contained in an I-cell from $Cov.$
In the paper we will consider mainly bounded horizontal and vertical
strips in ${\bf R}^2,$ their interval coverings and subdivisions. \\
Let $ {\bf X} = ({\bf x}_{1}, \cdots,{\bf x}_{n}) =
([{\underline x}_{1}, {\overline x}_{1}], \cdots,
[{\underline x}_{n}, {\overline x}_{n}] $ be the n-dimensional
real interval vector with
$ {\underline x}_{i} \leq x_{i} \leq {\overline x}_{i} $
 ("rectangle" or "box").
Let $f$ be a real continuous function  of  $n$ variables that is
defined on ${\bf X}, \; Of-$ its optimal interval
evaluation~\cite{AH:IC} on ${\bf X}$. The pair $({\bf X}, Of)$ is
called the {\em interval functional element.} If $Ef$ is an
interval that contains $Of$ then we will call the pair $({\bf
X},Ef)$ {\em the extension} of $({\bf X}, Of)$ or  $eif-${\em
element.}      \\ Let $f$ be the constant signs function on ${\bf
X}.$ If $f > 0$ (respectively $f < 0$) on ${\bf X}$ and $Of > 0$
(respectively $Of < 0$) then we will call $({\bf X}, Of)$ {\em the
correct interval functional element} (shortly $c-${\em element}).
More generally we will call {\em the correct interval functional
element} an extension $({\bf X},Ef)$ of $({\bf X}, Of)$ that has
the same sign as $Of.$


\section{Some Categories and Functors of Interval Mathematics}

A set of intervals with inclusion relation forms a category ${\cal CIP}$ of
preorder~\cite{Gl:IA}.
\begin{Def}
  A contravariant functor from ${\cal CIP}$ to the category of sets
is called the interval presheaf.
\end{Def}
  For a finite set $FS = \{{\bf X}_{i}\}$ of $m-$dimensional intervals in
${\bf R}^n, \; m \le n,$ the union $V$ of the intervals forms a
piecewise-linear
manifold in ${\bf R}^n.$ Let $G$ be the graph of the adjacency relation
of intervals from $FS$. The manifold $V$ is connected if $G$ is a
connected graph. In the paper we are considering connected manifolds.
Let $f$ be a constant signs function on
${\bf X} \in FS.$ The set $\{({\bf X}_{j},Of)\}$ of $c-$elements (if exists)
is called {\em a constant signs continuation of} $f$ on $\{{\bf X}_{j}\}.$ If
$\{{\bf X}_{j}\}$ is the maximal subset of $FS$ relatively a constant signs
function $f$ then $\{({\bf X}_{j},Of)\}$ is called the {\em constant signs
continuation of} $f$ on $FS$.

\section{On Interval Operads}

Operads was introduced  by J. May~\cite{M:GI}.
Operadic language is useful for investigation of many problems in
mathematics and physics. Some recent applications
of operads in physics is given by J. Morava~\cite{Mo:BT}.
 Here we give a short description of interval operad.
 The space of continuous interval
functions of $j$ variables forms the topological space $I{\cal C}(j).$
Its points are operations $I^{j}{\bf R} \to I{\bf R}$ of arity  $j.$
$I{\cal C}(0)$ is a single point $*.$ The class of interval spaces
$I^{n}{\bf R}, \; n \ge 0$ forms the category~\cite{Gl:IA}. We will consider
the spaces with base points and denote the category of those spaces by
$I{\cal U}.$  Let $X \in I{\cal U}$ and for $k \ge 0$ let $I{\cal E}(k)$
be the space of maps $M(X^{k},X).$ There is the action (by permuting the
inputs) of the symmetric group $S_k$ on $I{\cal E}(k).$ The identity
element $1 \in I{\cal E}(1)$ is the identity map of $X.$
\begin{PropDef}
 In the above mentioned conventions let $k \ge 0$ and $j_{1}, \ldots ,
j_{k} \ge 0$ be integers. Let for each choice of $k$ and  $j_{1}, \ldots ,
j_{k}$  there is a map
\begin{eqnarray*}
\gamma: I{\cal E}(k) \times I{\cal E}(j_{1}) \ldots \times I{\cal E}(j_{k})
\to I{\cal E}(j_{1} + \ldots + j_{k})
\end{eqnarray*}
given by multivariable composition. If maps $\gamma$ satisfy associativity,
equivalence and unital properties then $I{\cal E}$ is the endomorphism
interval operad $I{\cal E}_X$ of $X.$

\end{PropDef}

\section{Interval Programs as Functors}

Let $A = I{\bf R}$ be the interval algebra~\cite{AH:IC} with
interval arithmetic operations. In many cases extra interval operations
are required. So we have to extend the notion of interval algebra.
Let us define the "operator domain" $\Omega$ of interval computations
as sequence of sets $\Omega_0$ (interval constants and variables),
$\Omega_1$ (unary interval operations), $\Omega_2$
(binary interval operations). $\ldots$. In these notations the set
$T_{\Omega}$ of all "non-branching" programs in $\Omega$ is defined
as the least subset of $(\bigcup_{n=0}^{\infty} \Omega_{n} \bigcup
\{( )\})^*$ such that following axioms are satisfied: \\
(t) $\Omega_{0} \subseteq T_{\Omega};$     \\
(tt) for $n \ge 1, \; \omega \in \Omega_{n}$ and $t_{1}, \ldots ,
t_{n} \in T_{\Omega}, \; \omega(t_{1}, \ldots ,t_{n}) \in T_{\Omega}.$ \\
$\Omega-$interval algebra is constructed from the interval algebra $A$
and functions $\omega_{A}: \; A^{n} \to A, \; \omega \in \Omega_{n}.$
Below in the section our results follows Gougen~\cite{Go:CT} who discussed
the non-interval case.
\begin{Prop}
$T_{\Omega}$ is an initial object in the category of $\Omega-$
interval algebras.
\end{Prop}
 Let $IP$ be the interval program that implements an interval computation,
$G = G(IP)$ the graph of the flow diagram of $IP,$ ${\cal G}^{\otimes}$
the category
of all paths in $G.$ Let ${\cal IPF}$ be the category of interval sets with
partial interval functions.
\begin{Prop}
 Interval program $IP$ defines a functor $\overline{IP}: {\cal G}^{\otimes}
\to {\cal IPF}.$
\end{Prop}

\smallskip

{\bf \Large Some Problems}

\smallskip

Let $E$ be an interval expression with a tree $T.$
Let $w(E)$ be the diameter of $E$ on a given interval data.
\begin{Probl}
How to transform $E$ to the interval expression of the minimal complexity?
\end{Probl}

\begin{Probl}
How to transform $E$ to an expression with $min \; w(E)?$
\end{Probl}

Let $S$ be a sufficiently smooth bounded surface with boundary over
$xy-$plane in ${\bf R}^3.$ \\
Let $S_{xy} = Proj_{xy}S$ be the projection of $S$ on the $xy-$plane,
 $S_{x} = Proj_{x}S_{xy}, \; S_{y} = Proj_{y}S_{xy}$.
\begin{Probl}
Compute the set of points that minimize for each fixed $x$ the distance

$$ z_{min} = \{\min z(x,y) \; |  \; x \in S_{x}, y \in S_{y} \}.$$
\end{Probl}

I know a solution of the problem 3 for Minkowski-Cohn moduli space $(1)$
and for some another surfaces.

\subsection*{Acknowledgements}

I would like to thanks the organizers of the  NATO ASI -
AXIOMATIC, ENRICHED \& MOTIVIC HOMOTOPY THEORY (Isaac Newton
Institute (Cambridge)) and the  SCAN2002  for providing a very
pleasant environment during the conferences and respectively for
support and for partial support.



\end{document}